\documentclass[12pt]{article}
\usepackage[all]{xy}
\usepackage{amsmath}
\usepackage{amsfonts}
\usepackage{amssymb}
\usepackage{amscd}
\usepackage{amsthm}
\usepackage{latexsym}
\usepackage{amsbsy}

\newtheorem{lem}{Lemma}[section]
\newtheorem{prop}{Proposition}[section]
\newtheorem{theorem}{Theorem}[section]
\newtheorem{corollary}{Corollary}[section]
\renewcommand{\hom}{\mbox{\textnormal{Hom}}{\hspace{3pt}}}

\begin{document}
\date{}
\title
{Equivariant Cyclic Cohomology of H-Algebras}
\author{R. Akbarpour and  M. Khalkhali \\Department of Mathematics, 
University of Western Ontario\\
\texttt{akbarpur@uwo.ca \; masoud@uwo.ca}}
\maketitle

\begin{abstract}
We define an equivariant $K_0$-theory for \textit{Yetter-Drinfeld} algebras over 
a Hopf algebra with an invertible antipode. We then show that this definition can 
be generalized to all Hopf-module algebras. 
We show that there exists a pairing, generalizing Connes' pairing, between
this theory and a suitably defined Hopf algebra equivariant cyclic cohomology theory.

\end{abstract}


\section{Introduction}

Equivariant cyclic cohomology, for actions of discrete groups or compact Lie groups on algebras, has been studied by various authors
~\cite{bg94,b87,gj93,fs86,n89,n93}. One of the main themes studied in these papers is the relation between equivariant cyclic cohomology and the cyclic
cohomology of the corresponding crossed product algebra. In~\cite{rm01}, we extended one of the main results of these investigations, 
namely the Feigin-Tsygan and (independently) Nistor spectral sequence~\cite{fs86,n89} to actions of Hopf algebras. The $\mathsf{E_2}$-term
of this spectral sequence can be considered as the complex of noncommutative equivariant de Rham cochains on the given algebra.

One of the main features of cyclic cohomology is the existence of a pairing, through Connes' Chern character, between 
$K$-theory
and cyclic cohomology. In attempting to extend this pairing to a Hopf algebra equivariant setting one faces the 
following problem.
Let $\mathcal{H}$ be a Hopf algebra. For an $\mathcal{H}$-module algebra $A,$ here called an $\mathcal{H}$-algebra, and a finite dimensional $\mathcal{H}$-module
$V,$ one would like a natural algebra structure on $A \otimes End(V)$ together with a natural $\mathcal{H}$-action, to turn it into an $\mathcal{H}$-algebra.
There seems to be two possible approaches to this problem, depending on whether we use the diagonal $\mathcal{H}$-action or not. 
In Section 6 we show that the two definitions in fact coincide. If we use the diagonal $\mathcal{H}$-action on $A \otimes 
End(V)$, then the tensor product algebra structure on $A \otimes End(V)$ is not an $\mathcal{H}$-algebra, unless $\mathcal{H}$ 
is cocommutative. In general, one has to twist this tensor product structure with the help of an extra structure on $A$. This problem is
naturally solved by introducing the class of Yetter-Drinfeld 
algebras over a Hopf algebra. One can then show that for a Yetter-Drinfeld algebra $A$, $A \otimes End(V)$ has a natural
$\mathcal{H}$-algebra structure. In fact we first discovered formula~(\ref{eq:syd}) in Lemma~\ref{lem:yd} and realized 
later on that this 
condition is equivalent to a well known condition, namely the Yetter-Drinfeld condition, 
provided the antipode of $\mathcal{H}$
is bijective.
We define the equivariant $K_0$-theory, $K_0^{\mathcal{H}}(A)$, of a 
Yetter-Drinfeld algebra $A$ as the Grothendieck group
of the semigroup of invariant idempotents in $A \otimes End(V)$ for all finite dimensional $\mathcal{H}$-modules $V$. 

Alternatively, one can define a non-diagonal $\mathcal{H}$-action on $A \otimes End(V)$ by embedding it into 
$End(A \otimes V)$
first and then induce the 
conjugation $\mathcal{H}$-action as in~\cite{nt03} (see formula~(\ref{eq:A2})). One can then check that endowed with the tensor product algebra structure, 
$A \otimes End(V)$ is an $\mathcal{H}$-algebra. This approach leads to an apparently different definition of equivariant 
$K$-theory~\cite{nt03} that also admits a pairing with the equivariant cyclic cohomology defined in this paper.
In Section 6, motivated by examples from (co)quasitriangular Hopf algebras, we show that the two definitions in fact coincide.
We should also mention that the trace map used in~\cite{nt03} is exactly the trace map introduced in an earlier version of
the present paper for cocommutative Hopf algebras.
An interesting feature in our generalization of Connes' Chern character which is patterned after Connes' 
original construction in~\cite{aC85}  is the equivariant trace map $\Psi$ (Proposition~\ref{prop:psi}). 
We remark that the complex of cyclic equivariant cochains introduced in Section $3$ is not quite the same as the complex
that naturally appeared in~\cite{rm01}. We can prove, however, that it enjoys the same relation to crossed 
product algebras (Theorem~\ref{th:di}). This complex behaves better with respect to pairing with $K$-theory and 
this motivated our choice.

\section{Preliminaries}
 In this paper we work over a fixed field $k$ of characteristic zero. We denote the coproduct, antipode and counit of a 
Hopf algebra by $\Delta$, $S$ and
$\epsilon$, respectively.
Let $\mathcal{H}$ be a Hopf algebra. We use Sweedler's notation and write $\Delta h=h^{(0)} \otimes h^{(1)}$, where summation is understood.
Similarly, we write $\Delta^{(n)} h=h^{(0)} \otimes h^{(1)} \otimes \dots \otimes h^{(n)},$ where 
$\Delta^{(n)}:\mathcal{H} \rightarrow 
\mathcal{H}^{\otimes (n+1)}$ is defined by $\Delta^{(1)}=\Delta$ and $\Delta^{(n)}=
(\Delta \otimes 1) \circ \Delta^{(n-1)}$, $n \ge 2$.
By a left $\mathcal{H}$-module we mean a left $\mathcal{H}$-module over the underlying
algebra of $\mathcal{H}$. Let $A$ be an algebra. We say $A$ is a left $\mathcal{H}$-algebra, also called a left 
$\mathcal{H}$-module algebra, if $A$ is a left $\mathcal{H}$-
module and for all $a,b \in A, h \in \mathcal{H},$
\begin{eqnarray*}
& &h \cdot (ab)=(h^{(0)} \cdot a)(h^{(1)} \cdot b),\\
& &h \cdot 1=\epsilon(h)1.
\end{eqnarray*}

By a \textit{paracocyclic} object in a category $\mathcal{A}$~\cite{fs86,gj93} we mean a cosimplicial object $A$ in $\mathcal{A}$ endowed with 
operators $\tau_n:A_n \rightarrow A_n$, called cyclic operators, such that the following extra relations are satisfied: 
\begin{eqnarray}  
\tau_{n+1} \partial^i = \partial^{i-1} \tau_{n},\hspace{15pt} 1\le i \le n\; , \quad
\tau_{n+1} \partial^0  = \partial^{n+1},  \\
\tau_{n-1} \sigma^i  =  \sigma^{i-1} \tau_{n},  \hspace{15pt} 1\le i\le n\; , \quad
\tau_{n-1} \sigma^0  = \sigma^{n-1} \tau_{n}^2,  \notag
\end{eqnarray}
where $\partial^i:A_n \rightarrow A_{n+1}$ are coface maps and $\sigma^i:A_n \rightarrow A_{n-1}$ are codegenerecies.
We note that a similar notion is independently introduced in~\cite{n89,n93}. If in addition we have $\tau_n^{n+1}=id$ for all $n \ge 0$, then we have a \textit{cocyclic} object in the 
sense of Connes~\cite{aC994}.
By a \textit{bi-paracocyclic} object in $\mathcal{A}$, we mean a paracocyclic object in
the category of paracocyclic objects in $\mathcal{A}$. So, giving a bi-paracocyclic object
in  $\mathcal{A}$ is equivalent to giving a double sequence $A(p,q)$ of objects of $\mathcal{A}$ and operators
$\partial_{p,q} ,\sigma_{p,q},\tau_{p,q}$ and
$ \bar{\partial}_{p,q} ,\bar{\sigma}_{p,q},\bar{\tau}_{p,q}$ such that, for all $p \ge 0$,
\begin{eqnarray*}
B_p(q)=\{ A(p,q),\sigma^i_{p,q} ,\partial^i_{p,q} ,\tau_{p,q} \},
\end{eqnarray*}
and for all $q \ge 0$,
\begin{eqnarray*}
\bar{B}_q(p)=\{ A(p,q),\bar{\sigma}^i_{p,q} ,\bar{\partial}^i_{p,q} ,\bar{\tau}_{p,q} \},
\end{eqnarray*}
are paracocyclic objects in $\mathcal{A}$ and every horizontal operator commutes with every vertical operator.

We say that a bi-paracocyclic object is \textit{cocylindrical}~\cite{gj93} if for all $p,q \ge 0,$
\begin{eqnarray}  
\bar{\tau}_{p,q}^{p+1} \; \tau_{p,q}^{q+1}=id_{p,q}.
\end{eqnarray}

If $A$ is a bi-paracocyclic object in  $\mathcal{A}$, the paracocyclic
object related to the diagonal of $A$ will be denoted by $\Delta A$. So, the paracocyclic operators on $\Delta A(n)=A(n,n)$ are $\bar{\partial}^i_{n,n+1}
 \partial^i_{n,n},$ $\bar{\sigma}^i_{n,n-1} 
\sigma^i_{n,n},$ $\bar{\tau}_{n,n}  \tau_{n,n}.$
When $A$ is cocylindrical, since the cyclic operator of $\Delta A $ is $\bar{\tau}_{n,n} \tau_{n,n}$ and  $\bar{\tau},\tau$
commute, then, from
$\bar{\tau}^{n+1}_{n,n} \tau^{n+1}_{n,n} =id_{n,n}$, we conclude that $(\bar{\tau}_{n,n} \tau_{n,n})^{n+1}=id$. So that $\Delta A$ is a cocyclic object.

A \textit{paracochain complex}~\cite{gj93}, by definition, is a graded $k$-module $\mathsf{V}^{\bullet} = ( V^i)_{i \ge 0}$\; equipped
with operators $b: V^i \rightarrow V^{i+1}$ and $ B: V^i \rightarrow V^{i-1} $ such that $ b^2 = B^2 =0 ,$ and the operator
$T = 1-(bB+Bb)$ is invertible. In the case that $T=1$, the paracochain complex is called a \textit{mixed complex}.

Corresponding to any paracocyclic module $A$, we can define the paracochain complex $\mathsf{C}^{\bullet}(A)$ with the underlying graded 
module   $\mathsf{C}^n(A)= A(n)$ and the operators $ b = \sum_{i=0}^n (-1)^i \partial^i$ and 
$ B =  N \sigma (1 - (-1)^{n+1} \tau ) $. Here, $\sigma$ is the extra degeneracy satisfying $\tau \sigma^0 = \sigma \tau $, and
 $ N= \sum_{i=0}^n(-1)^{in} \tau^i $ is the norm operator. For any bi-paracocyclic module $A$, $Tot(\mathsf{C}(A))$
is a paracochain complex with $Tot^n(\mathsf{C}(A))= \sum_{p+q=n} A(p,q)$ and with the operators $ Tot(b) = b + \bar{b} $ and
$Tot(B)=  B + T \bar{B}$, where $ T = 1-(bB+Bb)$. It is a mixed complex if $A$ is cocylindrical~\cite{gj93}.

If we define the normalized cochain functor $\mathsf{N}$ from  paracocyclic modules to paracochain complexes with the underlying
graded module $\mathsf{N}^n(A) = {\bigcap}_{i=0}^{n-1} ker(\sigma^i) $ and the operators $b ,B$ induced from $\mathsf{C}^{\bullet}(A)$, then we have the 
following well-known results (see~\cite{gj93} for a dual version):\\
1. The inclusion $ ( \mathsf{N}^{\bullet}(A) , b) \rightarrow (\mathsf{C}^{\bullet}(A),b)$ is a quasi-isomorphism of complexes.\\
2. The cyclic Eilenberg-Zilber theorem holds for cocylindrical modules, i.e., 
for any cocylindrical module $A$, there is a natural quasi-isomorphism $\mathbf{f}_0 + \mathbf{u f_1} : \mathsf{N}^{\bullet}(\Delta(A))
 \rightarrow Tot^{\bullet}(\mathsf{N}(A))$ of mixed complexes, where $\mathbf{f}_0$ is the shuffle map.

\section{Equivariant cyclic cohomology of H-algebras} 
In this section we introduce the complex of cyclic equivariant cochains for $\mathcal{H}$-algebras. It is a 
noncommutative analogue of the complex of equivariant differential forms (Cartan model, see e.g.~\cite{bg94}). Since compact quantum
groups naturally coact on interesting algebras like quantum spheres, it would be perhaps more natural to consider Hopf comodule
algebras. Passing to this dual setting does not present serious difficulties. Our cocyclic module in Theorem~\ref{th:cy} is not quite the 
dual of the cyclic module that appeared in the $\mathsf{E}^2$-term of the spectral sequence in~\cite{rm01}, 
but is very similar to it. In particular, Theorem~\ref{th:di} in the next section shows that this version of equivariant cyclic 
cohomology enjoys the same relation with cyclic cohomology of crossed product algebras as in the main theorem of~\cite{rm01}.
The reason we prefer the present complex is that it works 
better for pairing with $K$-theory.
 
Let $\mathcal{H}$ be a Hopf algebra with a bijective antipode
and let $F(\mathcal{H})$ be the space of k-linear maps
$f: \mathcal{H}\to k$. Let $A$ be an $\mathcal{H}$-algebra and let $C^n(A,F(\mathcal{H}))$
denote the linear space of $(n+1)$-linear mappings
$$f: A^{\otimes (n+1)} = A\otimes A\otimes \cdots \otimes A \to F(\mathcal{H}).$$
We define an $\mathcal{H}$-action on $C^{n}(A,F(\mathcal{H}))$ by 
\begin{equation}\label{eq:ac}
(h\cdot f)(a_0 ,a_1 ,\dots ,a_n)(g) =f(h^{(0)}\cdot a_0 ,h^{(1)}\cdot a_1 ,\dots ,h^{(n)} \cdot a_n)(g) \hspace{15pt} h,g\in \mathcal{H},\; a_i \in A
\end{equation}

A cochain $f\in C^n (A,F(\mathcal{H}))$ is called \textit{$\mathcal{H}$-equivariant} if for all $h,g \in \mathcal{H}, a_i \in A,$
\begin{eqnarray}
& &h\cdot f(a_0 ,\dots ,a_n)(g) = f(a_0 ,\dots ,a_n)(S^{-1}(h)\cdot g) \\
&=& f(a_0 ,\dots ,a_n)(S(h^{(1)})gh^{(0)}), \notag
\end{eqnarray}
where the action on the left is defined in~(\ref{eq:ac}) and the action on the right is defined by
$$h\cdot g=S^2(h^{(0)})gS(h^{(1)}),$$ so that $S^{-1}(h)\cdot g=S(h^{(1)})g h^{(0)}.$ We define $C_\mathcal{H}^n(A)$ to be
the linear space of all $\mathcal{H}$-equivariant $f\in C^{n}(A,F(\mathcal{H}))$.

We define a cocyclic module structure on the spaces $\{C_\mathcal{H}^n(A)\}_{n \ge 0}$. First
we define the cyclic operator $T_n$ on $C_\mathcal{H}^n(A)$ by
\begin{equation}
T_nf (a_0 ,a_1 ,\dots ,a_n)(g) := f(S^{-1}(g^{(0)})\cdot a_n,a_0 ,\dots ,a_{n-1})
(g^{(1)}).
\end{equation}

\begin{lem}
$T_n $ is an equivariant map, i.e., $T_n f \in C_\mathcal{H}^n (A)$, for $f \in C_\mathcal{H}^n (A)$. 
\end{lem}
\begin{proof}
We should check that $h\cdot T_n f(g) = T_n f(S^{-1}(h)\cdot g)$:\\

\noindent $T_n f(a_0 ,a_1 ,\dots ,a_n)(S^{-1}(h) \cdot g)$
\begin{eqnarray*}
&=& T_n f(a_0 ,a_1 ,\dots ,a_n)
(S(h^{(1)})gh^{(0)} )\\
&=& f((S^{-1}(h^{(0)} )S^{-1}(g^{(0)}) h^{(3)} ) \cdot a_n ,a_0 ,\dots ,a_{n-1} )
(S(h^{(2)} )g^{(1)} h^{(1)} ) \\
&=& f((S^{-1}(h^{(0)} )S^{-1}(g^{(0)} )h^{(2)} ) \cdot a_n ,a_0 ,\dots ,a_{n-1} )
(S^{-1}(h^{(1)} ) \cdot g^{(1)} ) \\
&=& h^{(1)} \cdot f((S^{-1}(h^{(0)} )S^{-1}(g^{(0)} )h^{(2)} )\cdot a_n ,a_0 ,\dots ,
a_{n-1} )(g^{(1)} ) \\
&=& f((h^{(1)} S^{-1}(h^{(0)})S^{-1}(g^{(0)} )h^{(n+2)} )\cdot a_n ,h_{(2)}\cdot a_0 ,
\dots ,h^{(n+1)}\cdot a_{(n-1)})(g^{(1)}) \\
&=& f(S^{-1}(g^{(0)} )\cdot (h^{(n)} \cdot a_n) ,h^{(0)} \cdot a_0,\dots ,h^{(n-1)}
\cdot a_{n-1} )(g^{(1)} ) \\
&=& T_n f(h^{(0)} \cdot a_0 ,\dots ,h^{(n)} \cdot a_n)(g)\\
&=& h\cdot T_n f(a_0 ,\dots ,a_n)(g).
\end{eqnarray*}
\end{proof}
We define the coface and codegeneracy
operators on $C_\mathcal{H}^n (A)$ as follows:
$$\partial ^i:C_\mathcal{H}^{n-1} (A) \to C_\mathcal{H}^n (A), \;\;\;\;\; \sigma^i :C_\mathcal{H}^{n+1} (A)\to C_\mathcal{H}^n (A),$$
\begin{eqnarray} \label{eq:cy1}
& &\partial^i f(a_0 ,\dots ,a_n )(g) = f(a_0 ,\dots ,a_ia_{i+1} ,\dots ,a_n)(g),
\quad 0 \le i \le n-1,  \notag\\ 
& &\partial ^n f(a_0 ,\dots ,a_n)(g) = f((S^{-1}(g^{(0)} )\cdot a_n) a_0 ,a_1, \dots ,
a_{n-1} )(g^{(1)} ), \\ 
& &\sigma^i f(a_0 ,\dots ,a_{n-1} ,a_n)(g)=f(a_0 ,\dots ,a_i ,1,a_{i+1} ,\dots ,
a_n)(g), \;\; 0\le i \le n. \notag
\end{eqnarray}
One can check that these operators are well defined, i.e., they send equivariant cochains to equivariant cochains.
Now we are ready to state the main result of this section.
\begin{theorem} \label{th:cy}
For any Hopf algebra $\mathcal{H}$ with a bijective antipode
and an $\mathcal{H}$-algebra $A$, the $\mathcal{H}$-equivariant space
$C_\mathcal{H}^{\natural} (A) = \{ C_\mathcal{H}^n (A) \}_{n\ge 0}$ with operators defined in~(\ref{eq:cy1}),
is a cocyclic module.
\end{theorem}
\begin{proof}
We only check those identities that involve the cyclic operator $T$ and leave the rest to the reader.

\item[$\bullet \quad $] $T_n \partial^0 = \partial^n.$
\begin{eqnarray*}
& &(T_n \partial^0 f)(a_0 ,a_1 ,\dots ,a_n)(g)
= \partial^0 f(S^{-1}(g^{(0)} )\cdot a_n ,a_0 ,\dots ,a_{n-1} )(g^{(1)}) \\
&=& f((S^{-1}(g^{(0)} )\cdot a_n )a_0 ,\dots ,a_{n-1} )(g^{(1)} ) 
= \partial^n f(a_0 ,a_1 ,\dots ,a_n)(g).
\end{eqnarray*}
\item[$\bullet \quad $] $T_n \partial^i = \partial^{i-1} T_{n+1}$. For $\;1\le i < n$, this is obvious. For $i=n,$ we have
\begin{eqnarray*}
& &(T_n \partial^n f)(a_0 ,a_1 ,\dots ,a_n)(g)= 
\partial^n f(S^{-1}(g^{(0)} )\cdot a_n ,a_0 ,\dots ,a_{n-1} )(g^{(1)} ) \\
&=&
f((S^{-1}(g^{(1)} )\cdot a_{n-1} )(S^{-1}(g^{(0)}) \cdot a_n) ,a_0 ,\dots ,a_{n-2} )
(g^{(2)} )\\  
&=&
f(S^{-1}(g^{(0)} )\cdot (a_{n-1} a_n ),a_0 ,\dots ,a_{n-2} )(g^{(1)} ) \\
&=&
T_n f(a^0 ,\dots ,a^{n-2} ,a^{n-1} a^n )(g) 
=
(\partial^{n-1} T_{n+1} f)(a_0 ,\dots ,a_n)(g).
\end{eqnarray*}
\item[$\bullet \quad $] $T_n \sigma^0 = \sigma^n T_{n+1}^2.$
\begin{eqnarray*}
& &(\sigma^nT_{n+1}^2 f)(a_0 ,\dots ,a_n)(g)\\
&=& (T_{n+1}^2 f)(a_0 ,\dots ,a_n ,1)(g) 
= T_{n+1} f((S^{-1}(g^{(0)} )\cdot 1),a_0 ,\dots ,a_n)(g^{(1)} ), \\
&=& f((S^{-1}(g^{(1)} )\cdot a_n ),(S^{-1}(g^{(0)} )\cdot 1) , a_0 ,\dots ,a_{n-1} )
(g^{(2)} ), \\
&=& f(S^{-1}(g^{(0)}) \cdot a_n ,1,a_0 ,\dots ,a_{n-1} )(g^{(1)}) 
= (\sigma^0 f)(S^{-1}(g^{(0)} )\cdot a_n ,a_0 ,\dots ,a_{n-1} )(g^{(1)} ) \\
&=& (T_n \sigma^0 f)(a_0 ,a_1 ,\dots ,a_{n-1} ,a_n)(g).
\end{eqnarray*}

\item[$\bullet \quad $]
$T_n\sigma^i = \sigma^{i-1} T_{n+1}, \; 1\le i \le n.$ This is obvious.

\item[$\bullet \quad $] $T_n^{n+1} = \mbox{id}_n.$ Let $f\in C_\mathcal{H}^n (A)$. Then, 
\begin{eqnarray*}
& &T_n f(a_0 ,a_1 ,\dots ,a_n )(g) = f(S^{-1}(g^{(0)} ) \cdot a_n ,\dots ,a_{n-1} )
(g^{(1)} ) \\
& &T^2_n f(a_0 ,a_1 ,\dots ,a_n )(g)
= f(S^{-1}(g^{(1)} )\cdot a_{n-1} ,S^{-1}(g^{(0)} ) \cdot a_n ,\dots ,a_{n-2} )(g^{(2)} )\\
& &\vdots\\
& &T_n^n f(a_0 ,a_1 ,\dots ,a_n )(g)
= f(S^{-1}(g^{(n-1)}) \cdot a_1 ,S^{-1}(g^{(n-2} )
\cdot a_2 ,\dots ,S^{-1}(g^{(0)} ) \cdot a_n ,a_0 )(g^{(n)} ). 
\end{eqnarray*}
Thus,
\begin{eqnarray*}
& &T_n^{n+1} f(a_0 ,a_1 ,\dots ,a_n )(g)\\ 
&=& f(S^{-1}(g^{(n)} )\cdot a_0 ,S^{-1}(g^{(n-1)})
\cdot a_1 ,\dots ,S^{-1}(g^{(0)}) \cdot a_n )(g^{(n+1)} ) \\
&=& S^{-1}(g^{(0)} ) \cdot f(a_0 ,a_1 ,\dots ,a_n )(g^{(1)} ) 
= f( a_0 ,a_1 ,\dots ,a_n )(S^{-2}(g^{(0)}) \cdot g^{(1)} )\\ 
&=& f(a_0 ,a_1 ,\dots ,a_n )(g^{(0)} g^{(2)} S^{-1}(g^{(1)} )) 
= f(a_0 ,a_1 ,\dots ,a_n )(g^{(0)} \epsilon (g^{(1)} )) \\
&=& f( a_0 ,a_1 ,\dots ,a_n )(g).
\end{eqnarray*}

\noindent This last identity completes our proof of the theorem.
\end{proof}
We denote the Hochschild, cyclic, and periodic cyclic cohomology groups of the cocyclic module 
$\{C^n_{\mathcal{H}}(A)\}_{n \ge 0}$ by $HH_{\mathcal{H}}^{\bullet}(A)$, $HC^{\bullet}_{\mathcal{H}}(A)$, and
$HP_{\mathcal{H}}^{\bullet}(A)$, respectively.

\textbf{Example 1.}(trivial actions) Assume that $\mathcal{H}$ acts trivially on $A$, i.e., $h \cdot a= \epsilon(h) a$
for all $h \in \mathcal{H}, a \in A$. Then the cocyclic module $\{C^n_{\mathcal{H}}(A) \}_{n \ge 0}$ simplifies as
follows. We have $C^n_{\mathcal{H}}(A) \simeq C^n(A) \otimes R(\mathcal{H}),$ where $R(\mathcal{H}) \subset F(\mathcal{H})$
is the space of invariant linear functionals on $\mathcal{H}$. By definition, $f \in R(\mathcal{H})$ if $f(S^{-1}(h) 
\cdot g)= \epsilon(h) f(g)$ for all $h, g \in \mathcal{H}$. It then follows that $HC^n_{\mathcal{H}}(A) \simeq HC^n(A)
\otimes R(\mathcal{H})$, $n \ge 0$.

\textbf{Example 2.}(Morita invariance) Let $A$ be a left $\mathcal{H}$-algebra. Then the algebra of $r \times r$
matrices over $A$, $M_r(A)=A \otimes M_r(k)$, is a left $\mathcal{H}$-algebra where the left $\mathcal{H}$-action is
defined by $h \cdot (a \otimes m)=h \cdot a \otimes m$. The \textit{equivariant trace map} $tr:C^n_{\mathcal{H}}(A) 
\rightarrow C^n_{\mathcal{H}}(M_r(A)), n \ge 0,$ is defined by
$$
(trf)(a_0 \otimes m_0,\dots,a_n \otimes m_n)(g)= tr(m_0 \dots m_n)f(a_0, \dots,a_n)(g),
$$
where $tr$ is the usual trace on $M_r(k)$. It can be checked that $tr$ is a morphism of cocyclic modules. It follows from 
Corollary~\ref{mu}
that the induced map $tr:HC^n_{\mathcal{H}}(A) \rightarrow HC^n_{\mathcal{H}}(M_r(A))$ is an isomorphism for $n \ge 0$.

\textbf{Example 3.}
Let $\theta:A \rightarrow A$ be an automorphism of an algebra $A$. Then $A$ is an $k[x,x^{-1}]$-module algebra, where
$k[x,x^{-1}]$ is the Hopf algebra of Laurent polynomials. We identify the equivariant cyclic complex as follows. We have an 
isomorphism
\[
Hom(k[x,x^{-1}] \otimes A^{\otimes(n+1)},k) \simeq \prod_{-\infty}^{+\infty} Hom(A^{\otimes(n+1)},k), 
\]
sending $f \mapsto (f_m)^{+\infty}_{-\infty},$ where $f_m(a_0,\dots,a_n)=f(x^m,a_0,\dots,a_n).$ It is 
clear that a cochain $(f_m)^{+\infty}_{-\infty}$ is equivariant iff for all $m$
\[
f_m(\theta a_0,\theta a_1, \dots,\theta a_n)=f_m(a_0,\dots,a_n)\;\; \forall a_i \in A. 
\]
Thus we obtain a decomposition of cocyclic modules
\[
C_{\mathcal{H}}^n(A) \simeq \prod_{-\infty}^{+\infty} C_{\theta ,m}^n(A),
\]
where $C^n_{\theta,m}(A)=\{ f:A^{\otimes(n+1)} \rightarrow k, f(\theta a_0, \dots, \theta a_n)=f(a_0,\dots,a_n) \}.$
The coface and cyclic operators of $\{ C^n_{\theta,m}(A) \}_{n \ge 0}$ are given by
\begin{eqnarray*}
& &(\delta_i f)(a_0, \dots,a_n)=f(a_0, \dots,a_ia_{i+1}, \dots,a_{n+1}), \;\; 0\le i \le n,\\
& &(\delta_{n+1} f)(a_0, \dots,a_n)=f((\theta^m a_{n+1}) a_0,a_1, \dots,a_n),\\
& &(tf)(a_0, \dots,a_n)=f(\theta^m a_n, a_0, \dots, a_{n-1}).
\end{eqnarray*}
For $m=1$, the cocyclic module $\{C^n_{\theta,1}(A) \}_{n \ge 0}$ is exactly the cocyclic module in~\cite{kmt03}, 
used to define the  
"$\theta$-twisted cyclic cohomology" of $A$.

\section{Connection with the cocyclic module
$
\mathbf{Hom}_{ \mathbf{k} } ( ( \mathbf{A} \rtimes \mathbf{\mathcal{H}} ) ^ {\natural} ,\mathbf{k}) 
$
}
In this section we define a cocyclic map between the cocyclic module $C_\mathcal{H}^\natural (A)$ of equivariant 
cochains on $A$, introduced in
the previous section  and the cocyclic module $\hom_k ((A\rtimes
\mathcal{H})^{\natural}, k)$, associated with the crossed product algebra $A\rtimes \mathcal{H}$.                                              

Define a $k$-linear map 
$$\varphi_n : C_\mathcal{H}^n (A) \to \hom_k ((A\rtimes \mathcal{H})^{\otimes (n+1)} ,k), $$
by 

$\varphi_n f(a_0 \otimes g_0 ,\dots ,a_n\otimes g_n ):=$
\begin{multline*}
f(S^{-1}(g_0^{(0)} g_1^{(1)} g_2^{(2)} \cdots g_n^{(n)} )\cdot a_0 ,S^{-1}(g_1^{(0)}
g_2^{(1)} \cdots g_n^{(n-1)} )\cdot a_1 ,\\
\dots ,S^{-1}(g_{n-1}^{(0)} g_n^{(1)} )
\cdot a_{n-1} ,S^{-1}(g_n^{(0)} )\cdot a_n )(g_0^{(1)} g_1^{(2)} \cdots g_n^{(n)}).
\end{multline*}
Let $\varphi=\{\varphi_n\}_{n \ge 0}$.
Now we can state our first main result in this section.

\begin{theorem}
$\varphi$ defines a cocyclic map between cocyclic modules $C_\mathcal{H}^\natural (A)$
and $\hom_k ((A\rtimes \mathcal{H})^{\natural} ,k)$.
\end{theorem}

\begin{proof}

First we show that $\varphi$ commutes with cyclic operators. We have\\

$(\varphi_n T_n f)(a_0 \otimes g_0 ,\dots ,a_n\otimes g_n )$=
\begin{multline*}
T_n
f(S^{-1}(g_0^{(0)} g_1^{(1)} \cdots g_n^{(n)} )\cdot a_0 ,S^{-1}(g_1^{(0)}g_2^{(1)}\cdots
g_n^{(n-1)} )\cdot a_1 ,\\
\dots ,
S^{-1}(g_{n-1}^{(0)} g_n^{(1)} )\cdot a_{n-1} ,
S^{-1}(g_n^{(0)} )\cdot a_n )(g_0^{(1)} g_1^{(2)} \cdots g_n^{(n+1)} )
\end{multline*}
\begin{multline*}
\text{=} f((S^{-1}(g_0^{(1)} g_1^{(2)} \cdots g_n^{(n+1)} )S^{-1}(g_n^{(0)} )) \cdot a_n ,
S^{-1}(g_0^{(0)} g_1^{(1)} \cdots g_n^{(n)} )\cdot a_0 ,\\
\dots ,
S^{-1}(g_{n-1}^{(0)}
g_n^{(1)} )\cdot a_{n-1} )$ $(g_0^{(2)} g_1^{(3)} \cdots g_n^{(n+2)} )
\end{multline*}
\begin{multline*}
\text{=} f((S^{-1}(g_n^{(n+1)} )S^{-1}(g_n^{(0)} g_0^{(1)} \cdots g_{n-1}^{(n)} ))\cdot a_n ,
(S^{-1}(g_n^{(n)}) S^{-1}(g_0^{(0)} g_1^{(1)} \cdots g_{n-1}^{(n-1)} )\cdot a_0 ,\\
\dots ,
(S^{-1}(g_n^{(1)}) S^{-1}(g_{n-1}^{(0)} )) \cdot a_{n-1} )$ $(g_0^{(2)} g_1^{(3)}\cdots
g_n^{(n+2)} )
\end{multline*}
\begin{multline*}
\text{=} S^{-1}(g_n^{(1)} )\cdot f(S^{-1}(g_n^{(0)} g_0^{(1)} \cdots g_{n-1}^{(n)} ) \cdot a_n ,
S^{-1}(g_0^{(0)} g_1^{(1)} \cdots g_{n-1}^{(n-1)} )\cdot a_0 ,\\
\dots ,
S^{-1}(g_{n-1}^{(0)} )\cdot a_{n-1} )(g_0^{(2)} g_1^{(3)}\cdots
g_n^{(n+2)} )
\end{multline*}
\begin{multline*}
\text{=}f(S^{-1}(g_n^{(0)} g_0^{(1)} \cdots g_{n-1}^{(n)} ) \cdot a_n ,
S^{-1}(g_0^{(0)} g_1^{(1)} \cdots g_{n-1}^{(n-1)} )\cdot a_0 ,\\
\dots ,
S^{-1}(g_{n-1}^{(0)} )\cdot a_{n-1} )(S^{-2}(g_n^{(1)}) \cdot (g_0^{(2)} g_1^{(3)}\cdots
g_n^{(n+2)} )).
\end{multline*}
Since
\begin{eqnarray*}
& &S^{-2}(g_n^{(1)}) \cdot (g_0^{(2)} g_1^{(3)} \cdots g_n^{(n+2)}) = g_n^{(1)}
g_0^{(2)} g_1^{(3)} \cdots g_{n-1}^{(n+1)} g_n^{(3)} S^{-1}(g_n^{(2)} ) \\
&=& \epsilon (g_n^{(2)} )g_n^{(1)} g_0^{(2)} g_1^{(3)} \cdots g_{n-1}^{(n+1)} 
= g_n^{(1)} g_0^{(2)} g_1^{(3)} \cdots g_{n-1}^{(n+1)},
\end{eqnarray*}
we obtain

$(\varphi_n T_n f)(a_0 \otimes g_0 ,\dots ,a_n\otimes g_n )$=
\begin{multline*}
\text{=} f(S^{-1}(g_n^{(0)} g_0^{(1)} g_1^{(2)} \cdots g_{n-1}^{(n)} ) \cdot a_n ,
S^{-1}(g_0^{(0)} g_1^{(1)} \cdots g_{n-1}^{(n-1)}) \cdot a_0 ,\\
S^{-1}(g_1^{(0)} g_2^{(1)} \cdots g_{n-1}^{(n-2)}) \cdot a_1 ,
\dots ,
 S^{-1}(g_{n-1}^{(0)})
\cdot a_{n-1} )$ $(g_n^{(1)} g_0^{(2)} \cdots g_{n-1}^{(n+1)} ) .
\end{multline*}

On the other hand

\noindent $(\tau_n \varphi_n f)(a_0\otimes g_0 ,\dots ,a_n\otimes g_n)$ =$
\varphi_n f(a_n\otimes g_n ,a_0 \otimes g_0 ,\dots ,a_{n-1} \otimes g_{n-1}
)$ 
\begin{multline*}
\text{=} f(S^{-1}(g_n^{(0)} g_0^{(1)} g_1^{(2)} \cdots g_{n-1}^{(n)} )\cdot a_n ,
S^{-1}(g_0^{(0)} g_1^{(1)} \cdots g_{n-1}^{(n-1)} )\cdot a_0 ,\\
\dots ,S^{-1}(g_{n-2}^{(0)}
g_{n-1}^{(1)} )\cdot a_{n-2} ,S^{-1}(g_{n-1}^{(0)} )\cdot a_{n-1})(g_n^{(1)}
g_0^{(2)} \cdots g_{n-1}^{(n+1)} ).
\end{multline*}
Thus, $\varphi$ commutes with cyclic operators.

Next we show that $\varphi$ commutes with coface operators, i.e., 
$\partial^i \varphi_{n-1}=\varphi_n \partial^i.$ We check this only for $i=n$ and leave the rest to the reader.
We have\\

$(\partial^n \varphi_{n-1} f)(a_0 \otimes g_0 ,\dots , a_n \otimes g_n)$\\
=$(\varphi_{n-1} f)((a_n \otimes g_n)(a_0 \otimes g_0 ) ,a_1 \otimes g_1 ,\dots ,a_{n-1})  \otimes g_{n-1})$
\\
=$(\varphi_{n-1} f)(a_n (g_n^{(0)}) \cdot a_0 )\otimes g_n^{(1)} g_0 , a_1 \otimes g_1, \dots$
$, a_{n-1} \otimes g_{n-1})$
\begin{multline*}
\text{=} f(S^{-1}(g_n^{(1)} g_0^{(0)} g_1^{(1)} \dots g_{n-1}^{(n-1)})\cdot (a_n(g_n^{(0)} \cdot a_0)) ,S^{-1}(g_1^{(0)} g_2^{(1)}
\dots g_{n-1}^{(n-2)}) \cdot a_1 ,\\
\dots , S^{-1}(g_{n-1}^{(0)}) \cdot a_{n-1})(
g_n^{(2)} g_0^{(1)} \dots g_{n-1}^{(n)}),
\end{multline*}
and

$( \varphi_n \partial^n f)(a_0 \otimes g_0 , \dots , a_n \otimes g_n )=$
\begin{multline*}
\partial^n f ( S^{-1}(g_0^{(0)} g_1^{(1)} \dots g_n^{(n)}) \cdot a_0 , S^{-1}(g_1^{(0)} g_2^{(1)} \dots g_n^{(n-1)}) \cdot a_1 ,\\
\dots , S^{-1}(g_{n-1}^{(0)} g_n^{(1)}) \cdot a_{n-1} , S^{-1}(g_n^{(0)}) \cdot a_n)(g_0^{(1)} g_1^{(2)} \dots g_n^{(n+1)})
\end{multline*}
\begin{multline*}
\text{=} f(((S^{-1}(g_0^{(1)} g_1^{(2)} \dots g_n^{(n+1)})S^{-1}(g_n^{(0)})) \cdot a_n)(S^{-1}(g_0^{(0)}g_1^{(1)} \dots g_n^{(n)}) \cdot a_0),\\
S^{-1}(g_1^{(0)}g_2^{(1)} \dots g_n^{(n-1)}) \cdot a_1,
\dots , S^{-1}(g_{n-1}^{(0)}
g_n^{(1)}) \cdot a_{n-1})(g_0^{(2)}g_1^{(3)} \dots g_n^{(n+2)})
\end{multline*}
\begin{multline*}
\text{=} f((S^{-1}(g_n^{(0)}g_0^{(1)}g_1{(2)} \dots g_n^{(n+1)}) \cdot a_n)(S^{-1}(g_0^{(0)}g_1^{(1)} \dots 
g_n^{(n)}) \cdot a_0),\\S^{-1}(g_1^{(0)}g_2^{(1)} \dots g_n^{(n-1)}) \cdot a_1 ,
\dots , S^{-1}(g_{n-1}^{(0)}
g_n^{(1)}) \cdot a_{n-1})(g_0^{(2)}g_1^{(3)} \dots g_n^{(n+2)})
\end{multline*}
\begin{multline*}
\text{=} f(S^{-1}(g_0^{(0)}g_1^{(1)} \dots g_n^{(n)}) \cdot (S^{-1}(g_n^{(0)}) \cdot a_n)a_0),S^{-1}(g_1^{(0)}g_2^{(1)} \dots g_n^{(n-1)}) \cdot a_1,\\
\dots , S^{-1}(g_{n-1}^{(0)}g_n^{(1)}) \cdot a_{n-1})(g_0^{(1)}g_1^{(2)} \dots g_n^{(n+1)})
\end{multline*}
\begin{multline*}
\text{=} f((S^{-1}(g_0^{(0)}g_1^{(1)} \dots g_{n-1}^{(n-1)}g_n^{(n+1)})S^{-1}(g_n^{(1)} ))(a_n(g_n^{(0)} \cdot a_0)),\\
S^{-1}(g_1^{(0)}g_2^{(1)} \dots g_{n-1}^{(n-2)}g_n^{(n)}) \cdot a_1,
\dots,S^{-1}(g_{n-1}^{(0)}g_n^{(2)})
\cdot a_{n-1} )(g_0^{(1)}g_1^{(2)} \dots g_{n-1}^{(n)}g_n^{(n+2)})
\end{multline*}
\begin{multline*}
\text{=} S^{-1}(g_n^{(2)}) \cdot f (S^{-1}(g_n^{(1)}g_0^{(0)}g_1^{(1)} \dots g_{n-1}^{(n-1)}) \cdot (a_n(g_n^{(0)}
\cdot a_0)),\\S^{-1}(g_1^{(0)}g_2^{(1)} \dots g_{n-1}^{(n-2)}) \cdot a_1,
\dots , S^{-1}(g_{n-1}^{(0)}) \cdot
a_{n-1})(g_0^{(1)}g_1^{(2)} \dots g_{(n-1)}^{(n)}g_n^{(3)})
\end{multline*}
\begin{multline*}
\text{=} f(S^{-1}(g_n^{(1)}g_0^{(0)}g_1^{(1)} \dots g_{n-1}^{(n-1)}) \cdot (a_n(g_n^{(0)}) \cdot a_0)),
S^{-1}(g_1^{(0)}g_2^{(1)} \dots g_{n-1}^{(n-2)}) \cdot a_1,\\
\dots,S^{-1}(g_{n-1}^{(0)}) \cdot a_{n-1})
(g_n^{(2)}g_0^{(1)}g_1^{(2)} \dots g_{n-1}^{(n)}g_n^{(4)}S^{-1}(g_n^{(3)})).
\end{multline*}
Since $g_n^{(4)}S^{-1}(g_n^{(3)})=\epsilon(g_n^{(3)})$ and $\epsilon (g_n^{(3)})g_n^{(2)}=g_n^{(2)},$ 
the result is\\

$( \varphi_n \partial^n f)(a_0 \otimes g_0 , \dots , a_n \otimes g_n )  $
\begin{multline*}
\text{=} f(S^{-1}(g_n^{(1)}g_0^{(0)}g_1^{(1)} \dots g_{n-1}^{(n-1)}) \cdot (a_n(g_n^{(0)} \cdot a_0)),
S^{-1}(g_1^{(0)}g_2^{(1)} \dots g_{n-1}^{(n-2)}) \cdot a_1,\\
\dots , S^{-1}(g_{n-1}^{(0)}) \cdot a_{n-1})
(g_n^{(2)}g_0^{(1)}g_1^{(2)} \dots g_{n-1}^{(n)}).
\end{multline*}
The proof of compatibility of $\varphi$ with codegeneracies is similar and we leave it to the reader. The theorem is proved.

\end{proof}

\begin{corollary}
$\varphi$ induces natural maps between Hochschild, cyclic and
periodic cyclic cohomologies of $C_\mathcal{H}^{\natural} (A)$ and $Hom_k((A \rtimes \mathcal{H})^{\natural},k)$:
\[
\begin{CD}
HH^{\ast}_{\mathcal{H}}(A) @>{\varphi_H}>>HH^{\ast}(A \rtimes \mathcal{H}),\\
HC^{\ast}_{\mathcal{H}}(A) @>{\varphi_C}>>HC^{\ast}(A \rtimes \mathcal{H}),\\
HP^{\ast}_{\mathcal{H}}(A) @>{\varphi_P}>>HP^{\ast}(A \rtimes \mathcal{H}).
\end{CD}
\]
\end{corollary}

\begin{theorem} \label{th:di}
There is a spectral sequence that converges to the cyclic cohomology of $A \rtimes \mathcal{H}$. The $\mathsf{E}_2$-term
of this spectral sequence is given by
$$ \mathsf{E}_2^{p,q} = H^p(\mathcal{H},C^q_{\mathcal{H}}(A)).$$
\end{theorem}
\begin{proof}
We construct a cocylindrical module $X=\{X_{p,q}\}_{p,q \ge 0}$ and show that the diagonal $\Delta(X)$ of $X$ is isomorphic
to the cocyclic module $Hom((A \rtimes \mathcal{H})^{\natural}, k)$. We can then apply the cyclic Eilenberg-Zilber theorem to
derive our spectral sequence.
Let $$ X_{p,q}=Hom(A^{\otimes(p+1)} \otimes \mathcal{H}^{\otimes(q+1)},k).$$
We define the horizontal and vertical cosimplicial and cyclic operators by
\begin{eqnarray*}
\tau_{p,q}f(a_0,\dots,a_p)(g_0,\dots,g_q)=f(S^{-1}(g_0^{(0)} \dots g_q^{(0)}) \cdot a_p,a_0,\dots,a_{p-1})
(g_0^{(1)},\dots,g_q^{(1)}),\\
\partial_{p,q}^if(a_0,\dots,a_p)(g_0,\dots,g_q)=f(a_0,\dots,a_i a_{i+1},\dots,a_p)(g_0,\dots,g_q),\;0 \le i \le p-1,\\
\partial_{p,q}^pf(a_0,\dots,a_p)(g_0,\dots,g_q)=f((S^{-1}(g_0^{(0)} \dots g_q^{(0)}) \cdot a_p)a_0,\dots,a_{p-1})
(g_0^{(1)},\dots,g_q^{(1)}),\\
\sigma_{p,q}^if(a_0,\dots,a_p)(g_0,\dots,g_q)=f(a_0,\dots,a_i,1, a_{i+1},\dots,a_p)(g_0,\dots,g_q),\;0 \le i \le p,
\end{eqnarray*}
\begin{eqnarray*}
& &\bar{\tau}_{p,q}f(a_0,\dots,a_p)(g_0,\dots,g_q)=f(g_q^{(0)} \cdot (a_0,\dots,a_{p}))
(g_q^{(1)},g_0,\dots,g_{q-1}),\\
& &\bar{\partial}_{p,q}^if(a_0,\dots,a_p)(g_0,\dots,g_q)=f(a_0,\dots,a_p)(g_0,\dots,g_i g_{i+1},\dots,g_q),\;0 \le i \le q-1,\\
& &\bar{\partial}^q_{p,q}f(a_0,\dots,a_p)(g_0,\dots,g_q)=f(g_q^{(0)} \cdot (a_0,\dots,a_{p}))
(g_q^{(1)}g_0,\dots,g_{q-1}),\\
& &\bar{\sigma}_{p,q}^if(a_0,\dots,a_p)(g_0,\dots,g_q)=f(a_0,\dots,a_p)(g_0,\dots,g_i,1, g_{i+1},\dots,g_q),\;0 \le i \le q.\\
\end{eqnarray*}
One can check that $\{X_{p,q}\}_{p,q \ge 0}$ is a cocylindrical module.
The proof is very long, but is totally similar to the proof of Theorem 3.1 in~\cite{rm01} and is left to the reader.

Next we show that the diagonal of $X$, $\Delta(X)$, is isomorphic with the cocyclic module 
$Hom((A \rtimes \mathcal{H})^{\natural},k)$. To this end, we define the maps $\varphi=\{\varphi_n\}_{n \ge 0}$
and $\psi=\{\psi_n\}_{n \ge 0}$ by
\[
\varphi_n: \Delta^n(X) \rightarrow Hom((A \rtimes \mathcal{H})^{\otimes (n+1)},k),\;
\psi_n: Hom((A \rtimes \mathcal{H})^{\otimes (n+1)},k)  \rightarrow \Delta^n(X),
\]
\begin{multline*}
\varphi_n f(a_0 \otimes g_0 ,\dots ,a_n\otimes g_n )=\\
f(S^{-1}(g_0^{(0)} g_1^{(1)} g_2^{(2)} \cdots g_n^{(n)} )\cdot a_0 ,S^{-1}(g_1^{(0)}
g_2^{(1)} \cdots g_n^{(n-1)} )\cdot a_1 ,\\
\dots ,S^{-1}(g_{n-1}^{(0)} g_n^{(1)} )
\cdot a_{n-1} ,S^{-1}(g_n^{(0)} )\cdot a_n )(g_0^{(1)} ,g_1^{(2)}, \cdots, g_n^{(n)}),
\end{multline*}
$\psi_nf(a_0 , \dots ,a_n)(g_0, \dots , g_n)=$
\begin{eqnarray*}
f((g_0^{(0)} \dots g_{n-1}^{(0)} g_n^{(0)} ) \cdot a_0 \otimes g_0^{(1)} ,( g_1^{(1)} \dots g_{n-1}^{(1)}g_n^{(1)}) \cdot a_1 \otimes g_1^{(2)},
\dots , g_n^{(n)} \cdot a_n \otimes g_n^{(n+1)}).
\end{eqnarray*}

By a rather long computation one can verify that $\phi$ is a morphism of cocyclic modules 
and $\phi \circ \psi=\psi \circ \phi=id$.
Now, we can apply the generalized cyclic Eilenberg-Zilber theorem to derive our spectral sequence.
Again the argument is similar to that used in~\cite{rm01,gj93} and hence omitted.
\end{proof}

\begin{corollary} \label{co:sem}
Assume $\mathcal{H}$ is semisimple. Then we have an isomorphism of cyclic cohomology groups 
$HC^{\bullet}(A \rtimes \mathcal{H}) \simeq HC^{\bullet}_{\mathcal{H}}(A)$.
\end{corollary}
\textbf{Remark.} One can develop a similar theory for Hopf comodule algebras and prove the analogues of Theorem~\ref{th:di}
and Corollary~\ref{co:sem}, for cosemisimple Hopf algebras. It is known that compact quantum groups in the sense of 
Woronowicz are cosemisimple~\cite{ks97}.

We should also mention that some of our constructions and definitions in Sections 3 and 4, once appropriately dualized, reduce to those considered by Nistor in~\cite{n93}.
In particular, let $G$ be a compact Lie group acting smoothly on a complete locally convex algebra $A,$ and 
let $\mathcal{H}=Rep(G) \subset C^{\infty}(G)$ be the Hopf algebra of representable functions on $G$~\cite{ks97}. Then,
the dual of Corollary 4.2 reduces to Proposition 3.4 in~\cite{n93}.
\section{Equivariant K-Theory }

In this section we define the equivariant $K_0$-theory of Yetter-Drinfeld algebras and show that there exists a pairing, 
generalizing Connes' Chern character~\cite{aC85}, between this theory and the equivariant cyclic cohomology defined in 
Section 3. One can perhaps define an equivariant $K_0$-theory for any Hopf module algebra using finitely generated projective
modules endowed with a compatible action of the Hopf algebra. It is however not clear how to define a Chern character map
in this setting. Our approach, based on idempotents, however, naturally led us to a special class of $\mathcal{H}$-algebras,
namely the Yetter-Drinfeld $\mathcal{H}$-algebras.

Let $\mathcal{H}$ be a Hopf algebra with a bijective antipode. By a \textit{Yetter-Drinfeld} $\mathcal{H}$-algebra~\cite{ks97} we mean an 
algebra $A$ that satisfies the following conditions: 

1)$A$ is a left $\mathcal{H}$-algebra,

2)$A$ is a right $\mathcal{H}^{op}$-comodule algebra, i.e., the coaction $\rho:A \rightarrow A \otimes \mathcal{H},$ 
satisfies
\[\rho(ab)=a_{<0>} b_{<0>} \otimes b_{<1>} a_{<1>},\] where 
$\rho(a)=a_{<0>} \otimes a_{<1>} \in A \otimes \mathcal{H},$ denotes the coaction.

3) Conditions 1) and 2) are compatible in the sense that they satisfy the Yetter-Drinfeld condition
\begin{equation} \label{eq:com}
(h^{(1)} \cdot a)_{<0>} \otimes (h^{(1)} \cdot a)_{<1>}h^{(0)}= h^{(0)} \cdot a_{<0>} \otimes h^{(1)} a_{<1>},\; h \in \mathcal{H}, a\in A.
\end{equation}
We denote the 
class of Yetter-Drinfeld algebras of the above type by 
$ {_\mathcal{H}}\mathcal{YD}^{\mathcal{H}}$. It is easily checked that if $\mathcal{H}$ is cocommutative, then any 
left $\mathcal{H}$-algebra is a Yetter-Drinfeld algebra with a coaction defined by $a \rightarrow a \otimes 1$.

\begin{lem} \label{eq:l2}
Given any left $\mathcal{H}$-algebra $B$ and a Yetter-Drinfeld algebra $A$, then $A \otimes B$ with 
diagonal action and the following
multiplication is an $\mathcal{H}$-algebra:
\begin{equation} \label{eq:m1}
(a \otimes b)(c \otimes d)=ac_{<0>} \otimes (c_{<1>} \cdot b)d.
\end{equation} 
The $\mathcal{H}$-action on $A \otimes B$ is diagonal, i.e., 
$$h \cdot (a \otimes b)=h^{(0)} \cdot a \otimes h^{(1)} \cdot b.$$ 
\end{lem}
\begin{proof}
It is not difficult to see that~(\ref{eq:m1}) defines an associative product on $A \otimes B$. 
We check the $\mathcal{H}$-module algebra condition.
\begin{eqnarray*}
& &(h^{(0)} \cdot (a \otimes b))(h^{(1)} \cdot (c \otimes d))=
(h^{(0)} \cdot a \otimes h^{(1)} \cdot b)(h^{(2)} \cdot c \otimes h^{(3)} \cdot d)\\
&=&(h^{(0)} \cdot a)(h^{(2)} \cdot c)_{<0>} \otimes ((h^{(2)} \cdot c)_{<1>} \cdot 
( h^{(1)} \cdot b))(h^{(3)} \cdot d)\\
&=&(h^{(0)} \cdot a)(h^{(2)} \cdot c)_{<0>} \otimes ((h^{(2)} \cdot c)_{<1>}
h^{(1)}) \cdot b)(h^{(3)} \cdot d)\\
&=&(h^{(0)} \cdot a)(h^{(1)} \cdot c_{<0>}) \otimes ((h^{(2)} c_{<1>}) \cdot b) (h^{(3)} \cdot d)\\
&=&(h^{(0)} \cdot a)(h^{(1)} \cdot c_{<0>}) \otimes (h^{(2)}  \cdot (c_{<1>} \cdot b)) (h^{(3)} \cdot d)\\
&=&h \cdot( a c_{<0>} \otimes (c_{<1>} \cdot b) d)=h \cdot ((a \otimes b)(c \otimes d)). 
\end{eqnarray*}
\end{proof}

\begin{lem} \label{lem:yd}
Let $\mathcal{H}$ be a Hopf algebra with an invertible antipode and $A$ a Yetter-Drinfeld algebra. Then we have:
\begin{equation} \label{eq:syd}
\rho(h \cdot a)=h^{(1)} \cdot a_{<0>} \otimes h^{(2)} a_{<1>} S^{-1}(h^{(0)}).
\end{equation}
\end{lem}
\begin{proof}
By (\ref{eq:com}) we can see that
\begin{eqnarray*}
h^{(1)} \cdot a_{<0>} \otimes h^{(2)} a_{<1>} S^{-1}(h^{(0)}) &=&
(h^{(2)} \cdot a_{<1>})_{<0>} \otimes (h^{(2)} \cdot a_{<1>})_{<1>}h^{(1)} S^{-1}(h^{(0)})\\
&=&(h \cdot a)_{<0>} \otimes (h \cdot a)_{<1>}.
\end{eqnarray*}
\end{proof}
Conversely, one can check that condition~(\ref{eq:syd}) implies the Yetter-Drinfeld condition~(\ref{eq:com}).
  
Now let $V$ be a representation of $\mathcal{H}$, i.e., $V$ is a left $\mathcal{H}$-module with 
structure map $r:\mathcal{H} \rightarrow End(V)$. Then $B=End(V)$, with conjugation action 
$$ h \cdot f=r(h^{(0)}) \circ f \circ r(S(h^{(1)})),$$
is an $\mathcal{H}$-algebra.
Let $A \in {_\mathcal{H}}\mathcal{YD}^{\mathcal{H}}$. Then by Lemma~\ref{eq:l2}, $A \otimes End(V)$ is an $\mathcal{H}$- 
algebra with diagonal action and twisted multiplication i.e.,
\begin{eqnarray}\label{eq:A1}
& &(a \otimes u)(c \otimes v)=ac_{<0>} \otimes (c_{<1>} \cdot u)v,\;\;
h \cdot (a \otimes u)=h^{(0)} \cdot a \otimes h^{(1)} \cdot u. 
\end{eqnarray}
To simplify the notation, we denote the image of $h$ under $r$ by $h$ itself. 

Let $A$ be an $\mathcal{H}$-algebra. We say that $b \in A$ is an $\mathcal{H}$-invariant element if, for every $h \in \mathcal{H}$,
$h \cdot b = \epsilon (h) b.$ 
For a Yetter-Drinfeld algebra $A$ we define $P_\mathcal{H}(A)$ to be the set of all 
$\mathcal{H}$-invariant idempotents in all of the algebras $A \otimes End(V)$, where $V$ is a finite dimensional 
representation of $\mathcal{H}$. 
For $e , e' \in P_\mathcal{H}(A)$, $e \in A \otimes End(V)$ and  $e' \in A \otimes End(W)$, we define their sum 
$e_1 \oplus e_2$ as $\left(
\begin{smallmatrix}
e_1 & 0\\
0 & e_2
\end{smallmatrix}
\right) \in A \otimes End(V \oplus W)$.

Two $\mathcal{H}$-invariant idempotents $e \in A \otimes End(V)$ and $e' \in A \otimes End(W)$ are called Murray-von Neumann
equivalent if there exist $\mathcal{H}$-invariant elements $\gamma_1 \in A \otimes Hom_k(V,W) $ and 
$\gamma_2 \in A \otimes Hom_k(W,V) $ such that $\gamma_2 \gamma_1 = e$ and $\gamma_1 \gamma_2 = e'$.
Let $S_\mathcal{H}(A)$ denote the set of equivalance classes of $P_\mathcal{H}(A)$ under the Murray-von Neumann equivalence
relation. It is clear that $S_{\mathcal{H}}(A)$ is an abelian semigroup under the direct sum of idempotents. 

If there exists
an $\mathcal{H}$-invariant invertible element $\gamma \in Hom_k(V,W) \otimes A$ such that $\gamma e \gamma^{-1} = e'$, 
we say that $e$ and $e'$ are \textit{similar} and we write $ e \sim e' $.

The proof of the following theorem is similar to the case of group actions (Prop. 2.4.11 in~\cite{nC87}),
\begin{prop} \label{prop:sim}
$S_{\mathcal{H}}(A)$ is equal to the set of equivalence classes in $P_\mathcal{H}(A)$ for the equivalence relation generated by
similarity and the relation $p \sim p \oplus 0$.
\end{prop} 
We define the equivariant $K_0$-theory of a Yetter-Drinfeld algebra $A$ over a Hopf algebra $\mathcal{H},$ denoted by
$K_0^\mathcal{H}(A)$, as the Grothendieck group of $S_\mathcal{H}(A)$.

Let $(A^{\mathcal{H}})^{\times}$ be the group of invertible $\mathcal{H}$-invariant 
elements of $A$. 
Any element $b \in
 (A^{\mathcal{H}})^{\times}$ acts by conjugation on $C^n_\mathcal{H}(A)$ by the formula 
\begin{eqnarray*}
\alpha_b f (a_0,a_1,\dots,a_n)(g) = f(b a_0 b^{-1} ,b a_1 b^{-1},\dots ,b a_n b^{-1})(g),
\end{eqnarray*}
and any $\mathcal{H}$-invariant element $b \in A^{\mathcal{H}},$ acts by inner derivations by the formula
\begin{eqnarray*}
\delta_b f(a_0,a_1,\dots,a_n)(g) = \sum_{i \geq 0} f (a_0,\dots ,a_{i-1},[b,a_i] ,a_{i+1}, \dots ,a_n)(g).   
\end{eqnarray*}
\begin{lem}
$\alpha_b$ and $\delta_b$ are equivariant cocyclic maps.
\end{lem}
\begin{lem} \label{lem:zero}
$\alpha_b$ induces the identity map
and $\delta_b$ induces the zero map on $HH^*_{\mathcal{H}}(A),$ 
$HC^*_{\mathcal{H}}(A)$ and $HP^*_{\mathcal{H}}(A)$.   
\end{lem}
\begin{proof}
We see that the maps $\theta^i:C^{n+1}_\mathcal{H}(A) \to C^n_\mathcal{H}(A) , 0 \le i \le n$, where
\[
\theta^i f(a_0,a_1, \dots,a_n)(g) = f(a_0 b^{-1},ba_{1}b^{-1} ,\dots , ba_{i}b^{-1} ,b,a_{i+1}, \dots , a_n)(g), 
\]
define a presimplicial homotopy between $\theta^n \partial^{n+1} = \alpha_b$ and $\theta^0 \partial^0 = id$,
so that $\alpha_b$ is an isomorphism on Hochschild cohomology and hence on (periodic)cyclic cohomology.
In the case of $\delta_b$ we see that since 
\[
\partial^{n+1}f(a_0,a_1,\dots,a_n,b)(g)=f(S(g^{(0)}) \cdot b,a_0,\dots, a_n)(g^{(1)})=f(b,a_0,\dots,b_n)(g),
\]
the map 
\[
\varrho f(a_0,a_1, \dots,a_n)(g) = \sum_{ 0 \le i \le n} {(-1)}^i f (a_0,\dots ,a_{i}, b ,a_{i+1}, \dots ,a_n)(g), 
\]
defines the homotopy $b \varrho + \varrho b =- \delta_b$ between $\delta_b$ and $0$. So, $\delta_b$ is $0$ on
Hochschild cohomology and since $[b,1]=0$, in the normalized form $\varrho B + B \varrho = 0$ and 
$\varrho$ defines a homotopy between $\delta_b$ and $0$. Thus, $\delta_b$
is $0$ on cyclic and periodic cyclic cohomology. 
\end{proof}
Let $V$ be a finite dimensional $\mathcal{H}$-module. We construct the \textit{generalized trace map} 
$\Psi^n$ between $C^{n}_\mathcal{H} (A)$ and $C^{n}_\mathcal{H} (A \otimes End (V)) $. We define
$$\Psi^n: C^{n}_\mathcal{H} (A) \rightarrow C^{n}_\mathcal{H} (A \otimes End (V)),$$
\begin{eqnarray}\label{eq:tr} 
& &\Psi^n f(a_0 \otimes u_0, \dots , a_n\otimes u_n )(g)
=f(a_{0<0>} ,a_{1<0>}, \dots ,a_{n<0>})(g^{(1)}) \notag\\  
& &\hspace{5cm} tr(S(a_{0<1>})u_0 S(a_{1<1>}) u_1 \dots S(a_{n<1>}) u_n g^{(0)}). 
\end{eqnarray}
\begin{prop} \label{prop:psi}
$\Psi $ is an equivariant cocyclic map.
\end{prop}
\begin{proof} 
We first prove that $\Psi^n$ is equivariant, i.e., if $f$ is an equivariant cochain then
$h \cdot \Psi^n f(g) = \Psi^n f (S^{-1}(h) \cdot g) = \Psi^n f(S(h^{(1)})gh^{(0)})$: 
\begin{eqnarray*}
& &h \cdot (\Psi^n f)(a_0 \otimes u_0, \dots , a_n\otimes u_n )(g)\\
&=& \Psi^n f (h^{(0)} \cdot a_0 \otimes h^{(1)} \cdot u_0, \dots , h^{(2n)} \cdot a_n \otimes h^{(2n+1)}\cdot u_n )(g)\\
&=& f( h^{(1)} \cdot a_{0<0>},h^{(6)} \cdot a_{1<0>}, \dots ,h^{(5n+1)} \cdot a_{n<0>})(g^{(1)})\\
& &tr(S(h^{(2)}a_{0<1>}S^{-1}(h^{(0)}))(h^{(3)} u_0 S(h^{(4)}))S(h^{(7)}a_{1<1>}S^{-1}(h^{(5)}))(h^{(8)} u_1 S(h^{(9)})) \dots \\
& &\hspace{6.5cm} S(h^{(5n+2)}a_{n<1>}S^{-1}(h^{(5n)}))(h^{(5n+3)} u_n S(h^{(5n+4)}))g^{(0)})\\
&=& f( h^{(1)} \cdot a_{0<0>},h^{(2)} \cdot a_{1<0>}, \dots ,h^{(n+1)} \cdot a_{n<0>})(g^{(1)})\\
& &\hspace{4cm} tr(h^{(0)}S(a_{0<1>})u_0 S(a_{1<1>}) u_1 \dots S(a_{n<1>}) u_n S(h^{(n+2)}) g^{(0)})
\end{eqnarray*}
\begin{eqnarray*}
&=& h^{(1)} \cdot f( a_{0<0>}, a_{1<0>}, \dots ,a_{n<0>})(g^{(1)})\\
& &\hspace{4cm} tr(S(a_{0<1>})u_0 S(a_{1<1>}) u_1 \dots S(a_{n<1>}) u_n S(h^{(2)}) g^{(0)} h^{(0)})\\
&=& f( a_0, a_1, \dots , a_n)(S^{-1}(h^{(1)}) \cdot g)\\ 
& & \hspace{4cm}tr(S(a_{0<1>})u_0 S(a_{1<1>}) u_1 \dots S(a_{n<1>}) u_n S(h^{(2)}) g^{(0)} h^{(0)}) \\
&=& (\Psi^n f)(a_0 \otimes v_0, \dots , a_n\otimes v_n )( S(h^{(1)}) g h^{(0)} ).                                    
\end{eqnarray*}

Next we check that $\Psi$ is a cyclic map. First we check that $\Psi^n$ commutes with the cyclic operators:
\begin{eqnarray*}
& &( T_n \Psi^n f ) (a_0 \otimes u_0, a_1 \otimes u_1, \dots ,a_n \otimes u_n)(g)\\
&=& \Psi^n f  (S^{-1}(g^{(1)}) \cdot a_n \otimes S^{-1}(g^{(0)}) \cdot u_n, a_0 \otimes u_0, \dots ,a_{n-1} 
\otimes u_{n-1})(g^{(2)})
\end{eqnarray*}
\begin{eqnarray*}
&=&f(S^{-1}(g^{(3)}) \cdot a_n ,a_0 , \dots ,a_{n-1})(g^{(6)})\\
& &\hspace{2cm}tr(S(S^{-1}(g^{(2)}) a_{n<1>}S^{-2}(g^{(4)}))S^{-1}(g^{(1)})u_n g^{(0)}
S(a_{0<1>})u_0  \dots S(a_{n-1<1>}) u_{n-1} g^{(5)})\\
&=&f(S^{-1}(g^{(1)}) \cdot a_n ,a_0 , \dots ,a_{n-1})(g^{(2)})
tr(S(a_{0<1>})u_0 S(a_{1<1>}) u_1 \dots S(a_{n<1>}) u_n g^{(0)})\\
&=&(  \Psi^n T_n f ) (a_0 \otimes u_0, a_1 \otimes u_1, \dots ,a_n \otimes u_n)(g). 
\end{eqnarray*}
Now since $A$ is a right $\mathcal{H}^{op}$-comodule algebra, we have $\rho(ab)=a_{<0>}b_{<0>} \otimes b_{<1>}a_{<1>},$ and since 
\[
(a_i \otimes u_i)(a_{i+1} \otimes u_{i+1})= a_i a_{i+1<0>} \otimes (a_{i+1<1>} \cdot u_i)u_{i+1},
\]
we can see that
\begin{eqnarray*}
S((a_i a_{i+1<0>})_{<1>})((a_{i+1<1>} \cdot u_i)u_{i+1})&=&S( a_{i+1<1>}a_{i<1>} )(a_{i+1<2>} u_i S(a_{i+1<3>})u_{i+1})\\
&=&S(a_{i<1>} ) u_i S(a_{i+1<1>}) u_{i+1}.
\end{eqnarray*}
Therefore $\Psi^n$ commutes with coface operators $\partial^i, 0 \le i <n$. Since $T_n \partial^0=\partial^n,$ $\Psi^n$
commutes with all coface operators and it is easy to check that it also commutes with codegeneracy operators. This proves 
the proposition.
\end{proof}

\begin{corollary}(Morita invariance) \label{mu}
Let $V$ be a trivial $\mathcal{H}$-module. Then 
$HC^{\ast}_{\mathcal{H}}(A) \simeq HC^{\ast}_{\mathcal{H}}(A \otimes End(V))$.
\end{corollary}
\begin{proof}
As is shown in~\cite{mk94}, Theorem 6, Morita invariance is a formal consequence of two facts: inner automorphisms induce
the identity map on cohomology and a generalized trace map exists. In our case, these are established in Lemma 5.4 and Proposition 5.2.
\end{proof}

Now we state the main result of this section. We define $R(\mathcal{H})$ to be the space of
all invariant functions from $\mathcal{H}$ to the ground field k. So $f \in R(\mathcal{H})$ iff 
$f (h \cdot g) = f( S^2(h^{(0)}) g S(h^{(1)})) =\epsilon(h) f (g).$

\begin{theorem}
For each $n \ge 0$ there exists a bilinear pairing 
$K_0^\mathcal{H}(A) \times HC_{\mathcal{H}}^{2n}(A) \rightarrow R(\mathcal{H})$,
defined by
\begin{eqnarray*}
\langle [e] ,(f) \rangle (g) = \Psi^{2n} f(e, e , \dots , e)(g).
\end{eqnarray*}
We also have a pairing    
$K_0^{\mathcal{H}}(A) \times HP^{0}_\mathcal{H}(A) \rightarrow R(\mathcal{H})$, defined by
\begin{eqnarray*}
\langle [e] ,(f) \rangle (g) = \Psi f_0(e)(g)+\sum_{n} (-1)^n \frac{(2n)!}{n!} \Psi f_{2n}(e-\frac{1}{2}, e , \dots , e)(g),
\end{eqnarray*}
where in the last pairing $(f)=(f_{2n})_{0 \le n \le m}$ is an equivariant even periodic cyclic cocycle in 
the normalized equivariant $(b,B)$-bicomlex of $A$.

\end{theorem}

\begin{proof}
To check that the first pairing is well defined, let $f$ be a coboundary in $C_\mathcal{H}^{2n}(A)$. Then 
$\Psi^{2n} f = b \psi^{2n-1}$ is also a coboundary and we see that
\begin{eqnarray*}
& &\Psi^{2n} f(e,e, \dots ,e)(g) = b \psi^{2n-1}(e,e, \dots ,e)(g)
= \sum_{0 \le i \le n} (-1)^i \partial^i \psi^{2n-1} (e,e, \dots ,e)(g)\\
&=& \psi^{2n-1} (e^2,e, \dots ,e)(g)- \psi^{2n-1} (e,e^2, \dots ,e)(g)+ \dots +(-1)^{2n-2} \psi^{2n-1} 
(e,e, \dots ,e^2)(g) \\
& &\hspace{6cm}+(-1)^{2n-1} \psi^{2n-1} ((S^{-1}(g^{(0)}) \cdot e)e,e, \dots ,e)(g^{(1)})=0,
\end{eqnarray*}
since $e$ is $\mathcal{H}$-invariant and therefore $S^{-1}(g^{(0)}) \cdot e = \epsilon(g^{(0)}) e $.

Let $[e]=[e'] $ in $ K^\mathcal{H}_0(A)$, where $e \in A \otimes End(V)$ and $e' \in A \otimes End(W)$.
Then, by Proposition~\ref{prop:sim}, $e \sim e \oplus 0 \sim 0 \oplus e' \sim e'$,
so there exists an $\mathcal{H}$-invariant invertible element $\gamma \in A \otimes End(V \oplus W)$ such that
$e' \oplus 0 = \gamma ( 0 \oplus e ) \gamma^{-1}$. Then, by Lemma~\ref{lem:zero}, we have
\[
\Psi^{2n}f (e',\dots , e' )(g) = \Psi^{2n} f(\gamma e \gamma^{-1}, \dots ,\gamma e \gamma^{-1} )(g) = \Psi^{2n}f (e, \dots ,e)(g).
\]
Also, since f is equivariant, $\langle [e] ,(f) \rangle \in R(\mathcal{H})$. This 
finishes the proof of the first part. The proof of the second part is also similar to the non equivariant case 
as in~\cite{jL92} and is left to the reader.
\end{proof}


\section{Examples}
In this section we first give some examples of Yetter-Drinfeld module algebras and show that for Yetter-Drinfeld
module algebras naturally defined by $R$-matrix of (co)quasitriangular Hopf algebras the two definitions of 
$K_0^{\mathcal{H}}$ as given in this paper and in~\cite{nt03} coincide. We then generalize this result and show that the 
two definitions coincide for all Yetter-Drinfeld module algebras. Finally, we show that the quantum analogue of the Dirac
monopole line bundle over the quantum sphere $S^2_q$ defines an element of $U_q(su_2)$-equivariant $K$-theory of $S^2_q$.

Let $A$ be an $\mathcal{H}$-module algebra. Then one can check that $A \otimes \mathcal{H}$ with the following
structure is a Yetter-Drinfeld module algebra over $\mathcal{H}$:
\[
(a \otimes h)(b \otimes g)=(ab \otimes hg),
\]
\[ 
g \cdot(a \otimes h)= g^{(1)} \otimes g^{(0)} h S(g^{(2)}),\;\;\;\rho(a \otimes h) = 
a \otimes h^{(1)} \otimes S^{-1}(h^{(0)}).
\]
In particular $\mathcal{H}$ is a Yetter-Drinfeld module algebra over itself with the action and coaction defined as:
\[
g \cdot h= g^{(1)} h S(g^{(2)}),\;\;\;\rho(h) =h^{(1)} \otimes 
S^{-1}(h^{(0)}).
\]

Examples of Yetter-Drinfeld module algebras can be obtained also by considering  
(co)quasi-triangular Hopf algebras. It is
shown in~\cite{coz97} that given any $\mathcal{H}$-algebra (resp. $\mathcal{H}^{op}$-comodule algebra) $A$ over a quasitriangular
(resp. coquasitriangular) Hopf algebra $\mathcal{H},$ then $A$ can be turned into a Yetter-Drinfeld module algebra. We recall this
construction.

By definition, a \textit{quasitriangular Hopf algebra} is a pair $(\mathcal{H}, R)$, where $\mathcal{H}$ is a Hopf algebra and
$R=R^{(1)} \otimes R^{(2)}  \in \mathcal{H} \otimes \mathcal{H}$ is an invertible element which satisfies the following 
relations (R= r):
\[
\Delta(R^{(1)}) \otimes R^{(2)} = R^{(1)} \otimes r^{(1)} \otimes R^{(2)} r^{(2)},\;\;
 R^{(1)} \otimes  \Delta(R^{(2)}) = R^{(1)} r^{(1)} \otimes R^{(2)} \otimes r^{(2)},
\] 
\[
\Delta^{cop}(h)R=R \Delta(h),\;\; \epsilon(R^{(1)}) R^{(2)} =1,\;\;R^{(1)} \epsilon(R^{(2)}=1,
\]
\[
(S \otimes id)R=R^{-1},\;\; (id \otimes S)R^{-1}=R,\;\;(S \otimes S)R=R.
\]
for all $h \in \mathcal{H}$.

Given any left $\mathcal{H}$-module algebra $A,$ one can define a right $\mathcal{H}^{op}$-coaction on $A$ as follows:
\[
\rho(a)= R^{(2)} \cdot a \otimes R^{(1)}.
\]
It is easily checked that~\cite{coz97}, with the above coaction $A$ is a Yetter-Drinfeld $\mathcal{H}$-algebra.

A \textit{coquasitriangular Hopf algebra} is a pair $(\mathcal{H},\mathcal{R})$, where $\mathcal{H}$ is a Hopf algebra and 
 $\mathcal{R} \in (\mathcal{H} \otimes \mathcal{H})^{\ast}$ is a convolution-invertible map in the sense that there exists a map 
$\mathcal{R}^{-1} \in (\mathcal{H} \otimes \mathcal{H})^{\ast}$ such that
\[
\mathcal{R}^{-1}(h^{(0)} \otimes g^{(0)}) \mathcal{R}(h^{(1)} \otimes g^{(1)}) =
\mathcal{R}(h^{(0)} \otimes g^{(0)}) \mathcal{R}^{-1}(h^{(1)} \otimes g^{(1)})=
 \epsilon(h) \epsilon(g),
\]
and the following relations are satisfied:
\begin{eqnarray*}
& &\mathcal{R}(hg,r)=\mathcal{R}(h \otimes r^{(0)}) \mathcal{R}(g \otimes r^{(1)}),\;\;
\mathcal{R}(h,gr)=\mathcal{R}(h^{(0)} \otimes g) \mathcal{R}(h^{(1)} \otimes r),\\
& &g^{(0)} h^{(0)} \mathcal{R}(g^{(1)} \otimes h^{(1)})=\mathcal{R}(g^{(0)} \otimes  h^{(0)}) g^{(1)} h^{(1)},\;\;
\mathcal{R}(h \otimes 1)= \mathcal{R}(1 \otimes h)=\epsilon(a),\\
& &\mathcal{R}(S(h) \otimes g)=\mathcal{R}^{-1}(h \otimes g), 
\;\;\mathcal{R}^{-1}(h \otimes S(g))=\mathcal{R}(h \otimes g),\;\;
\mathcal{R}(S(h) \otimes S(g))=\mathcal{R}(h \otimes g),  
\end{eqnarray*}
for all $h,g,r \in \mathcal{H}.$

Now for any right $\mathcal{H}^{op}$-comodule algebra $A$ there is a left $\mathcal{H}$-module structure on $A$ defined by
\[
h \cdot a=a_{<0>} \mathcal{R}(h \otimes a_{<1>}),
\]
which turn it into a Yetter-Drinfeld module algebra.

Now let $\mathcal{H}$ be a quasitriangular Hopf algebra. The following lemma shows that the equivariant $K$-theory 
of the resulting Yetter-Drinfeld module algebra is independent of the choice of the $R$-matrix.
Recall from~\cite{nt03} that for any $\mathcal{H}$-module algebra $A$ and an $\mathcal{H}$-module $V$ there is an 
$\mathcal{H}$-module algebra structure on $A \otimes End(V)$ where the algebra structure is diagonal and the action 
is non-diagonal:
\begin{equation}\label{eq:A2}
(h \otimes u)(g \otimes v)=hg \otimes uv,\;\;h \cdot (a \otimes u)= h^{(1)} \cdot a \otimes h^{(0)} u S(h^{(2)}),
\end{equation}
for $a \in A$ and $u \in End(V)$. We denote this latter structure by $A \bar{\otimes} End(V)$ and our original $\mathcal{H}$
-algebra structure by $A \otimes End(V)$ as defined in~(\ref{eq:A1}).

\begin{lem}Let $A$ be an $\mathcal{H}$-algebra over a quasitriangular Hopf algebra and let $V$ be a left $\mathcal{H}$-module.
Then there is an $\mathcal{H}$-algebra isomorphism between $A \bar{\otimes} End(V)$  and $A \otimes End(V)$. 
\end{lem}
\begin{proof}
We prove that the following maps define $\mathcal{H}$-isomorphisms between $A \bar{\otimes} End(V)$ and 
$A \otimes End(V)$ inverse to each other:
\[
t:A \bar{\otimes} End(V) \rightarrow A \otimes End(V),\;t(a \otimes u)=R^{(2)} \cdot a \otimes R^{(1)} u,
\]
\[
t':A \otimes End(V) \rightarrow A \bar{\otimes} End(V),\;t'(a \otimes u)=R^{(2)} \cdot a \otimes S(R^{(1)}) u,
\]
where $R$ is the $R$-matrix of $\mathcal{H}$.
Since $t \circ t'= r^{(2)}R^{(2)} \cdot a \otimes S(r^{(1)})R^{(1)} u$ and 
$t' \circ t= R^{(2)}r^{(2)} \cdot a \otimes S(R^{(1)})r^{(1)} u$ and $(S \otimes id)R=R^{-1},$ therefore $t \circ t'=
t' \circ t=id$, where $R=r$.
Since 
\begin{eqnarray*}
& &t(a \otimes u)t(b \otimes w)=(R^{(2)} \cdot a \otimes R^{(1)} u)(r^{(2)} \cdot b \otimes r^{(1)} w)\\
& &=(R^{(2)} \cdot a)((r_1^{(2)} r^{(2)}) \cdot b) \otimes (r_1^{(1)} \cdot (R^{(1)} u)(r^{(1)}w)\\
& &=(R^{(2)} \cdot a)((r_1^{(2)} (\underbrace{r_2^{(2)} r^{(2)}})) \cdot b) \otimes r_1^{(1)} R^{(1)} u 
\underbrace{S(r_2^{(1)})r^{(1)}}w\\
& &=(R^{(2)} \cdot a)(r_1^{(2)} \cdot b) \otimes r_1^{(1)} R^{(1)} uw=R^{(2)} \cdot (ab) \otimes R^{(1)} uw=
t((a \otimes u)(b \otimes w)),
\end{eqnarray*}
where $R=r=r_1=r_2$, we see that $t$ is an algebra map. Also since
\begin{eqnarray*}
t(h \cdot (a \otimes u))=t(h^{(1)} \cdot a \otimes h^{(0)} u S(h^{(2)}))=
(R^{(2)}h^{(1)}) \cdot a \otimes R^{(1)}h^{(0)} u S(h^{(2)}),
\end{eqnarray*}
and $\Delta^{cop}(h)R=R \Delta(h)$, we conclude that
\begin{eqnarray*}
t(h \cdot (a \otimes u))=
(h^{(0)}R^{(2)}) \cdot a \otimes h^{(1)}R^{(1)} u S(h^{(2)})=
h^{(0)} \cdot (R^{(2)} a) \otimes h^{(1)} \cdot (R^{(1)} u)=h \cdot t(a \otimes u),
\end{eqnarray*}
which shows that $t$ preserves the $\mathcal{H}$-actions.
\end{proof}
A similar result holds in the coquasitriangular case. Motivated by these examples, we were led to the interesting
fact that our original $\mathcal{H}$-algebra structure on $A \otimes End(V)$ is always independent
of the choice of coaction:
\begin{prop}Let $A$ be an $\mathcal{H}$-Yetter-Drinfeld module algebra and $V$ be a representation of $\mathcal{H}$.
Then there is an $\mathcal{H}$-algebra isomorphism between $A \bar{\otimes} End(V)$with $\mathcal{H}$-algebra structure 
defined by~(\ref{eq:A2}) and $A \otimes End(V)$ with $\mathcal{H}$-algebra structure defined 
by~(\ref{eq:A1}).
\end{prop}
\begin{proof}
We prove that the following maps define an $\mathcal{H}$-isomorphism between $A \bar{\otimes} End(V)$ and 
$A \otimes End(V),$ inverse to each other:
\begin{eqnarray*}
\beta:A \bar{\otimes} End(V) \rightarrow A \otimes End(V),\;\beta(a \otimes u)=a_{<0>}  \otimes a_{<1>} u,\\
\beta':A \otimes End(V) \rightarrow A \bar{\otimes} End(V),\;\beta'(a \otimes u)=a_{<0>}  \otimes S(a_{<1>}) u.
\end{eqnarray*}
Since $\beta \circ \beta'( a \otimes u)=a_{<0>} \otimes a_{<1>} S(a_{<2>})u$ and 
$\beta' \circ \beta( a \otimes u)=a_{<0>} \otimes S(a_{<1>}) a_{<2>}u,$ we obtain $\beta \circ \beta'=\beta' \circ \beta=id$.
Now since
\begin{eqnarray*}
& &\beta(a \otimes u) \beta(b \otimes w)= a_{<0>} b_{<0>} \otimes (b_{<1>} \cdot (a_{<1>}u)(b_{<2>}w)\\
& &=a_{<0>} b_{<0>} \otimes b_{<1>} a_{<1>} u S(b_{<2>}) b_{<3>} w=a_{<0>}b_{<0>} \otimes a_{<1>}b_{<1>} uw=
\beta((a \otimes u)(b \otimes w)),
\end{eqnarray*}
$\beta$ is an algebra map. To check $\beta$ preserves the $\mathcal{H}$-actions, we see that
\begin{eqnarray*}
& &\beta(h \cdot (a \otimes u))=\beta(h^{(1)} \cdot a \otimes h^{(0)} u S(h^{(2)}))=
h^{(2)} \cdot a_{<0>} \otimes h^{(3)} a_{<1>} S^{-1}(h^{(1)}) h^{(0)} u S(h^{(4)})\\
& &=h^{(0)} \cdot a_{<0>} \otimes h^{(1)} a_{<0>} u S(h^{(2)})=h \cdot \beta(a \otimes u).
\end{eqnarray*}
\end{proof}
The above proposition shows that the two apparently different definitions of $K_0^{\mathcal{H}}$ in this paper and in~\cite{nt03}
are in fact the same. Moreover, under the isomorphism $\beta$, our generalized trace map~(\ref{eq:tr}) transforms to the 
following map
\begin{eqnarray}\label{eq:tr1} \notag 
& &(\Psi^n \circ \beta) f(a_0 \otimes u_0, \dots , a_n\otimes u_n )(g)
=\Psi^n  f(a_{0<0>} \otimes a_{0<1>}u_0, \dots , a_{n<0>} \otimes a_{n<1>}u_n )(g)\\ \notag 
& =&f(a_{0<0>} ,a_{1<0>}, \dots ,a_{n<0>})(g^{(1)})\\  
& &\hspace{2cm}tr(S(a_{0<1>})a_{0<2>}u_0 S(a_{1<1>})a_{1<2>} u_1 \dots S(a_{n<1>}) a_{n<2>}u_n g^{(0)})\\ \notag
&=&f(a_0 ,a_1, \dots ,a_n)(g^{(1)}) tr(u_0  u_1 \dots u_n g^{(0)}), 
\end{eqnarray}
which is the one used in~\cite{nt03}. In fact, the above trace map is exactly the same trace map that we introduced in the
first version of this paper for cocommutative Hopf algebras.

Now we show that the quantum monopole line bundle over the Podle\'s quantum sphere $S^2_q$ fits very well in our framework
and defines a $U_q(su_2)$-invariant idempotent. Recall that~\cite{p87} the Podle\'{s} equator quantum sphere $S^2_q$ is the 
$\ast$-algebra generated over $\mathbb{C}$ by the elements $a, a^{\ast}$ and $b$ subject to the relations
\[
a a^{\ast}+q^{-4}b^2=1,\;a^{\ast}a+b^2=1,\;ab=q^{-2}ba,\;a^{\ast}b=q^2ba^{\ast}.
\]

The quantum enveloping algebra $U_q(su_2)$ is a Hopf algebra over $\mathbb{C}$ generated by the 
elements $E,F,K$~\cite{ks97}
subject to the relations:
\[
KK^{-1}=K^{-1}K=1,\;\;KEK^{-1}=qE,\;\;KFK^{-1}=q^{-1}F,\;\;[F,E]=\frac{K^2-K^{-2}}{q-q^{-1}},
\]
\[
\Delta(K)=K \otimes K,\;\Delta(F)=F \otimes K+K^{-1} \otimes F,\;\Delta(E)=E \otimes K+K^{-1} \otimes E,
\]
\[
S(K)=K^{-1},\;S(E)=-qE,\;S(F)=-q^{-1}F,\; \epsilon(K)=1,\;\epsilon(E)=\epsilon(F)=0.
\]

By a direct computation one can show that:
\begin{lem}The following formulas define a $U_q(su_2)$-module algebra structure on $S^2_q$,
$$K \cdot a=qa,\;\;K \cdot a^{\ast}=q^{-1} a^{\ast},\;\;K \cdot b=b,$$
$$E \cdot b= q^{\frac{5}{2}} a,\;\; 
E \cdot {a^{\ast}} = - q^{\frac{3}{2}} (1+q^{-2}) b,\;\; E \cdot a=0,$$ 
$$F \cdot {a} = q^{-\frac{7}{2}} (1+q^2) b,\;\; F \cdot b=-q^{-\frac{1}{2}} a^{\ast},\;\; F \cdot a^{\ast}=0.$$
\end{lem}
The quantum analogue of the Dirac(or Hopf) monopole line bundle over $S^2$ is given by the following idempotent 
in $M_2(S^2_q)$
$$\mathbf{e}_q=\frac{1}{2}\left[\begin{array}{cc} 1+q^{-2}b & qa \\ q^{-1} a^{\ast} & 1-b \end{array}\right].$$

It can be directly checked that $\mathbf{e}_q^2=\mathbf{e}_q.$
Now we consider a 2-dimensional representation of $\mathcal{H}=U_q(su_2)$ on $V=\mathbb{C}^2$ defined by~\cite{ks97}
\[
E=\left[\begin{array}{cc} 0 & 0 \\ 1 & 0 \end{array}\right],\;F=\left[\begin{array}{cc} 0 & 1 \\ 0 & 0 \end{array}\right],\;
K=\left[\begin{array}{cc} q^{-\frac{1}{2}} & 0 \\ 0 & q^{\frac{1}{2}}  \end{array}\right].
\]
According to this representation $\mathbf{e}_q$ will be represented as
$$\mathbf{e}_q=\frac{1}{2}\left[\begin{array}{cc} 1+q^{-2}b & qa \\ q^{-1} a^{\ast} & 1-b \end{array}\right]=\frac{1}{2}
( 1 \otimes 1 + q^{-2} b \otimes FE - b \otimes EF + q a \otimes F + q^{-1} a^{\ast} \otimes E).$$

By a direct computation we show that $\mathbf{e}_q \in M_2(A)=A \bar{\otimes} End(V)$ is an $\mathcal{H}$-invariant idempotent, i.e., for all
$h \in U_q(su_2)$, $h \cdot \mathbf{e}_q=\epsilon(h) \mathbf{e}_q$. It suffices to check this for the generators $E,F$ and $K$, 
i.e., $E \cdot \mathbf{e}_q=0,F \cdot \mathbf{e}_q=0,K \cdot \mathbf{e}_q=\mathbf{e}_q.$

Since 
$$\Delta^2(K)=K \otimes K \otimes K,\;\;$$ 
$$\Delta^2(F)=F \otimes K \otimes K+K^{-1} \otimes F \otimes K+K^{-1} \otimes 
K^{-1} \otimes F,$$
$$\Delta^2(E)=E \otimes K \otimes K+K^{-1} \otimes E \otimes K+K^{-1} \otimes 
K^{-1} \otimes E,$$
we can see that
\begin{eqnarray*}
& &K \cdot \mathbf{e}_q=\frac{1}{2} [\; K \cdot 1 \otimes KK^{-1} + q^{-2} K \cdot b \otimes KFEK^{-1} - K \cdot b 
\otimes KEFK^{-1} + q \underbrace{K \cdot  a}_{qa} \otimes \underbrace{KFK^{-1}}_{q^{-1}F} \\
& &+ q^{-1} \underbrace{K \cdot a^{\ast}}_{q^{-1} a^{\ast}} \otimes \underbrace{KEK^{-1}}_{qE}\;]=\mathbf{e}_q,
\end{eqnarray*}
and
\begin{eqnarray*}
& &F \cdot \mathbf{e}_q = \frac{1}{2} [\; 
K \cdot 1 \otimes FK^{-1} + \underbrace{q^{-2} K \cdot b \otimes F^2 E K^{-1}}_{0}- K 
\cdot b \otimes FEFK^{-1} + \underbrace{q K \cdot a \otimes F^2 K^{-1}}_{0} +\\ 
& & q^{-1} K \cdot a^{\ast} \otimes FE K^{-1}
+\underbrace{F \cdot 1 \otimes K^{-2}}_{0} + q^{-2} F \cdot b \otimes K^{-1}F E K^{-1}- F
\cdot b \otimes K^{-1} EFK^{-1}+ \\
& & q F \cdot a \otimes K^{-1} F K^{-1} + \underbrace{q^{-1} F \cdot a^{\ast} \otimes K^{-1} E K^{-1}}_{0}
+K^{-1} \cdot 1 \otimes K^{-1} (-q^{-1}F)+\\
& &q^{-2} K^{-1} \cdot b \otimes K^{-1} FE (-q^{-1} F)- \underbrace{K^{-1} \cdot b \otimes K^{-1} EF (-q^{-1} F)}_{0} 
+ \underbrace{q K^{-1} \cdot a \otimes K^{-1} F (-q^{-1} F)}_{0}\\ 
& &+ q^{-1} K^{-1} \cdot a^{\ast} \otimes K^{-1} E (-q^{-1} F)\; ]. 
\end{eqnarray*}
Since
\begin{eqnarray*}
& &K \cdot 1 \otimes FK^{-1}=q^{-1}1 \otimes K^{-1}F=q^{-\frac{1}{2}} 1 \otimes F\\
& &K^{-1} \cdot 1 \otimes K^{-1} (-q^{-1}F)=-q^{-\frac{1}{2}} 1 \otimes F\\
& &-K \cdot b \otimes FEFK^{-1}=-q^{-1} b \otimes K^{-1}F=-q^{-\frac{1}{2}} b \otimes F\\
& &q^{-1} K \cdot a^{\ast} \otimes FE K^{-1}=q^{-2} a^{\ast} \otimes K^{-1} FE=q^{-\frac{3}{2}} a^{\ast} \otimes FE\\
& &q^{-2} F \cdot b \otimes K^{-1}F E K^{-1}=-q^{-\frac{5}{2}} a^{\ast} \otimes K^{-2} FE=-q^{-\frac{3}{2}} 
a^{\ast} \otimes FE\\
& &-F \cdot b \otimes K^{-1} EFK^{-1}=q^{-\frac{1}{2}} a^{\ast} \otimes K^{-2}EF=q^{-\frac{3}{2}} a^{\ast} \otimes EF,\\
& &q E \cdot a \otimes K^{-1} F K^{-1}=q^{-\frac{7}{2}}(1+q^2) b \otimes K^{-2} F= q^{-\frac{5}{2}}(1+q^2) b \otimes F\\
& &q^{-2} K^{-1} \cdot b \otimes K^{-1} FE (-q^{-1} F)=-q^{-3}b \otimes K^{-1} F=-q^{-\frac{5}{2}} b \otimes F,\\
& &q^{-1} K^{-1} \cdot a^{\ast} \otimes K^{-1} E (-q^{-1} F)=-q^{-1} a^{\ast} \otimes K^{-1}EF=-q^{-\frac{3}{2}}a^{\ast} 
\otimes EF,
\end{eqnarray*}
by adding the above relations, we obtain $F \cdot \mathbf{e}_q=0.$
Also we can see that
\begin{eqnarray*}
& &E \cdot \mathbf{e}_q = \frac{1}{2} [\; 
K \cdot 1 \otimes EK^{-1} + q^{-2} K \cdot b \otimes E F E K^{-1}- \underbrace{K 
\cdot b \otimes E^2 FK^{-1}}_{0} + q K \cdot a \otimes E F K^{-1}\\
& &+\underbrace{q^{-1} K \cdot a^{\ast} \otimes F^2 K^{-1}}_{0}
+\underbrace{E \cdot 1 \otimes K^{-2}}_{0} + q^{-2} E \cdot b \otimes K^{-1}F E K^{-1}- E
\cdot b \otimes K^{-1} EFK^{-1} +\\
& &\underbrace{q E \cdot a \otimes K^{-1} F K^{-1}}_{0} 
+ q^{-1} E \cdot a^{\ast} \otimes K^{-1} E K^{-1}
+K^{-1} \cdot 1 \otimes K^{-1} (-q E) +\\ 
& &\underbrace{q^{-2} K^{-1} \cdot b \otimes K^{-1} FE (-q E)}_{0}
-K^{-1} \cdot b \otimes K^{-1} EF (-q E) + q K^{-1} \cdot a \otimes K^{-1} F (-q E)+\\ 
& &\underbrace{q^{-1} K^{-1} \cdot a^{\ast} \otimes K^{-1} E (-q E)}_{0}\; ], 
\end{eqnarray*}
where,
\begin{eqnarray*}
& &K \cdot 1 \otimes EK^{-1}=q 1 \otimes K^{-1}E=q^{\frac{1}{2}} 1 \otimes E,\\
& &q^{-2} K \cdot b \otimes E F E K^{-1}=q^{-1} b \otimes K^{-1}E=q^{-\frac{3}{2}} b \otimes E,\\
& &q K \cdot a \otimes E F K^{-1}=q^2 a \otimes K^{-1} EF =q^{\frac{3}{2}} a \otimes EF,\\
& &q^{-2} E \cdot b \otimes K^{-1}F E K^{-1}=q^{\frac{1}{2}} a \otimes K^{-2}FE = q^{\frac{3}{2}} a \otimes FE,\\
& &-E\cdot b \otimes K^{-1} EFK^{-1}=-q^{\frac{5}{2}} a \otimes K^{-2}EF=-q^{\frac{3}{2}} a \otimes EF,\\
& &q^{-1} E \cdot a^{\ast} \otimes K^{-1} E K^{-1}=-q^{\frac{3}{2}}(1+q^{-2}) b \otimes K^{-2}E=
-q^{\frac{1}{2}}(1+q^{-2}) b \otimes E\\
& &K^{-1} \cdot 1 \otimes K^{-1} (-q E)-q 1\otimes K^{-1}E=-q^{\frac{1}{2}} 1 \otimes E,\\
& &-K^{-1} \cdot b \otimes K^{-1} EF (-q E)=q b \otimes K^{-1}E=q^{\frac{1}{2}} b \otimes E,\\
& &q K^{-1} \cdot a \otimes K^{-1} F (-q E)=-q a \otimes K^{-1}FE=-q^{\frac{3}{2}} a \otimes FE,
\end{eqnarray*} 
and therefore $E \cdot \mathbf{e}_q=0.$


\end{document}